\documentclass{amsart}
\usepackage{graphicx}

\newtheorem{thm}{Theorem}[section]

\newtheorem{prop}[thm]{Proposition}
\theoremstyle{definition}
\newtheorem{defn}[thm]{Definition}
\theoremstyle{remark}
\newtheorem{rem}[thm]{Remark}
\numberwithin{equation}{section}
\newcommand{\leftbul}[1]{\ensuremath{\bullet #1}}
\newcommand{\rightbul}[1]{\ensuremath{#1 \bullet}}
\newcommand{\dual}[1]{\ensuremath{\vec{#1}}}
\begin{document}

\title[meanders]{meanders in a Cayley graph}
\author{H. Tracy Hall}
\email{h.tracy@gmail.com}


\begin{abstract}
A meander of order $n$ is a simple closed curve in the plane
which intersects a horizontal line transversely at $2n$ points.
(Meanders which differ by an isotopy of the line and plane are
considered equivalent.)

Let $\Gamma_n$ be the Cayley graph of the symmetric group $S_n$ as
generated by all
$\binom{n}{2}$
transpositions.  Let $\Lambda_n$
be any interval of maximal length in $\Gamma_n$; this graph
is the Hasse diagram of the lattice of
noncrossing partitions.  The meanders of order
$n$ are in one-to-one correspondence with ordered pairs of maximally
separated vertices of $\Lambda_n$.

\end{abstract}
\maketitle

\begin{section}{introduction}
The lattice of noncrossing partitions and its variants have received
a fair amount of attention in recent years as a crossing grounds of
several seemingly unrelated disciplines of mathematics, including Artin
groups, cluster algebras, and free probability.  We here record a
connection between noncrossing partitions and the set of meanders,
or isotopy classes of simple closed curves
transversely intersecting an infinite line in the plane.
It is well recognized that meanders correspond to certain pairs
of noncrossing partitions, and indeed attempts have been made to
make use of this description to aid in the notoriously difficult
meander enumeration problem.  To our knowledge, however, it is a
new observation that the pairs of noncrossing partitions which give
rise to meanders are exactly those whose path distance is maximal
in the Hasse diagram of the lattice of noncrossing partitions---which
graph can itself be regarded as a geometric interval
in the Cayley graph of the symmetric group generated by all
two-element cycles.  It remains to be seen whether this
is any aid to the enumeration problem; the difficulty seems to stem
from the fact that the Hasse diagram is properly a directed graph,
and minimal paths connecting elements can ``change direction''
an arbitrary number of times.

In Section~2 we define meanders, noncrossing matchings and partitions,
and the Cayley graph $\Gamma_n$ and lattice $\Lambda_n$.  In
Section~3 we show that $\Lambda_n$ is a geometric interval in
$\Gamma_n$ and that meanders correspond to pairs of vertices
in $\Lambda_n$ whose path-distance is maximal, and finish
with some unanswered questions arising from this bijection.

The author wishes to thank Reinhard Franz, who first introduced
me to the problem of meander enumeration, as well as
organizers of combinatorics seminars at
Brigham Young University and the University of California, Berkeley,
where these results were first presented several years ago.  The
2005 Workshop "Braid Groups, Clusters and Free Probability" at
the American Institute of Mathematics was valuable not only
as a source of new information and new questions but also as
an encouragement to write up older answers
and make them known to the broader mathematical community.
Further encouragement came from David Savitt, who recently
observed and demonstrated independently \cite{S06} the main
result: that the number of components of a system
of meanders is measured exactly by path distances in the
lattice of noncrossing partitions.




\end{section}
\begin{section}{Preliminaries}
\begin{subsection}{Meanders}
A meander of order $n$ is a simple closed curve in the plane which
intersects an infinite horizontal line transversely in $2n$ points,
which points we label $\leftbul{1}, \rightbul{1}, \leftbul{2}, \rightbul{2},
\ldots, \leftbul{n}, \rightbul{n}$.
(Two meanders which differ only by isotopies of the line and the
upper and lower half-planes are considered equivalent, making
meanders combinatorial objects.)
For example, Figure~\ref{8meanders}
depicts the eight possible meanders of order~3.
We will refer to the collection of
points $\{\leftbul{1}, \leftbul{2}, \ldots, \leftbul{n} \}$ as {\em odd} points,
and the points $\{\rightbul{1}, \rightbul{2}, \ldots, \rightbul{n}\}$ we will call
{\em even}.
Meanders have a long history and still await a satisfactory enumeration;
see for example \cite{LZ93}.

\begin{figure} 
 \centering 
 \includegraphics[width=0.9\textwidth]{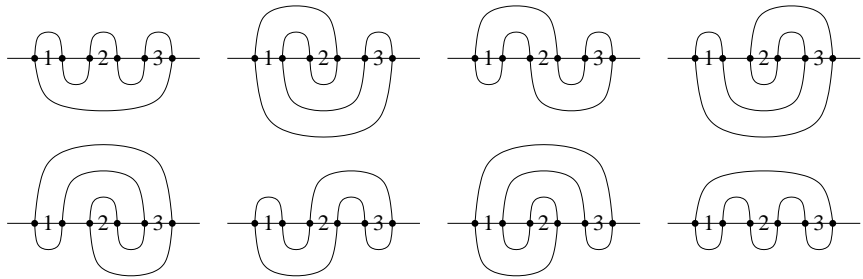} 
 \caption{The eight meanders of order $3$.} 
 \label{8meanders} 
\end{figure}

A meander is characterized by the labels of the endpoints
of arcs both above and below the dividing line, and since
every arc encloses either both endpoints or neither endpoint
of every other arc, each arc encloses an even total number
of points other than its own endpoints, and hence one of those
endpoints is even and one is odd.
There is thus a unique permutation
$\sigma$ of the integers $1, \ldots, n$ such that
$\leftbul{\sigma(i)}$
is an endpoint of an arc in the upper half-plane
if and only if the other endpoint of the arc is
$\rightbul{i}$,
%
%
and similarly a unique permutation
$\tau$ arising from the endpoints of arcs in the lower
half-plane.  The pair of permutations $(\sigma, \tau)$
is a complete invariant of meanders of a given order,
but not every pair of permutations arises
in this way.

Our principal concern is to characterize which
pairs of permutations $(\sigma, \tau)$ define a
meander, starting from topological considerations
in the plane (namely: no pair of arcs may intersect,
and the union of arcs must have a single component)
to arrive at a characterization
in terms of path lengths in
a Cayley graph of the symmetric group.

\end{subsection}

\begin{subsection}{Non-crossing partitions}
We first consider the collection of arcs on a single
side of the dividing line.
Any permutation $\sigma$ gives a matching
of even points $\rightbul{i}$ to odd points
$\leftbul{\sigma(i)}$, but for $\sigma$ to
represent half of a meander, there must exist
disjoint arcs in a single half-plane which
connect the pairs of matched points.
Such a system of $n$ arcs we call a
{\em noncrossing matching} of order $n$.  Every
such system
is canonically equivalent to a noncrossing partition,
which we define as follows:

\begin{defn}
Let $n$ points on the boundary of a circular disc
be labeled $1$ through $n$ in counter-clockwise order,
and let $\sim$ be an equivalence relation on the
set of marked points.  The partition arising
from the relation $\sim$ is said to be {\em noncrossing} if no
two blocks of the partition ``cross''---that is,
if whenever two points $i \sim k$ and two
other points $j \sim \ell$ have labels satisfying
$i < j < k < \ell$, then $j \sim k$ as well and
all four points are part of a single block.
Equivalently,
form the convex hull in the disc of each block
of points,
and the partition is noncrossing
if and only if the convex hulls are disjoint.
\end{defn}

\begin{figure} 
 \centering 
 \includegraphics[width=0.65\textwidth]{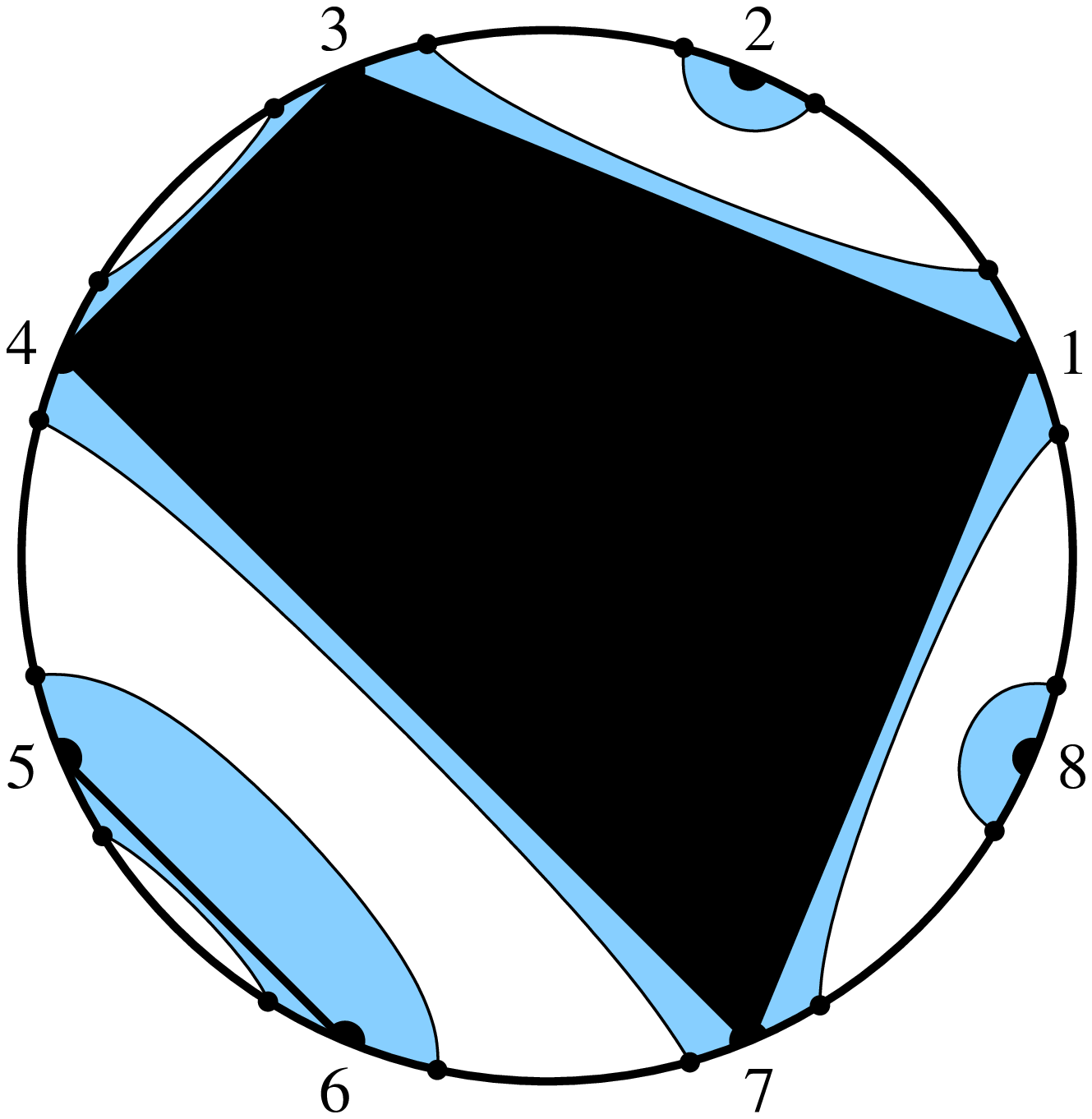} 
 \includegraphics{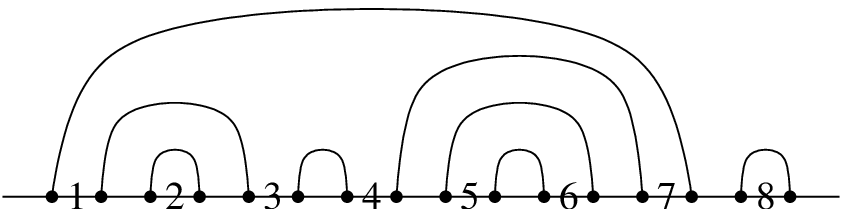}
 \caption{The noncrossing partition of order 8
coming from the permutation (1347)(56),
with convex hulls in black.
The $16$ points
$\leftbul{1}$, $\rightbul{1}$, etc.
are paired by the boundary arcs of the shaded region.} 
 \label{hulls} 
\end{figure}

\begin{prop} The noncrossing matchings of
order $n$
in a single half-plane
are in canonical bijection with noncrossing
partitions of order $n$.
\end{prop}

\begin{proof}
To begin with, we take the one-point compactification of the
half-plane (together with a reversal of orientation in the case
of the lower half-plane) so that the arcs join points
(marked $\leftbul{1}, \rightbul{1}, \ldots,
\leftbul{n}, \rightbul{n}$ in counter-clockwise order)
on the
boundary of a disc.
These $n$ arcs divide the disc into $n + 1$ connected components.
We define a relation $i \sim j$ if the section of disc boundary
which joins $\leftbul{i}$ to $\rightbul{i}$ is in the same connected
component as the section joining $\leftbul{j}$ to $\rightbul{j}$;
the relation $\sim$ defines a noncrossing partition of order $n$.
Conversely, if we start with a noncrossing partition of $n$ points
marked $1, \ldots, n$ in counter-clockwise order around a disc, form
convex hulls of the blocks of the partition, and take a regular
neighborhood of the union of convex hulls, the boundary of the
regular neighborhood
will consist of $n$ properly embedded
disjoint arcs in the disc, whose endpoints surround the points
$\{1, 2, \ldots, n\}$ in pairs
$(\leftbul{1}, \rightbul{1}), \ldots, (\leftbul{n}, \rightbul{n})$.
(Figure \ref{hulls} illustrates this correspondence for the example
of the noncrossing partitions whose blocks are $\{1,3,4,7\}$,
$\{2\}$, $\{5,6\}$, and $\{8\}$.)
\end{proof}


\end{subsection}

\begin{subsection}{The lattice $\Lambda_n$ of noncrossing partitions}
The set of all noncrossing partitions of a given order $n$ we
call $\Lambda_n$.  This set has a great deal of structure and
symmetry; in particular there is an operation
which reverses the parity of the points $\leftbul{1}, \rightbul{1}$,
etc., and which induces a symmetry of order $2n$ on $\Lambda_n$.
To every noncrossing partition we assign a dual as
follows:
\begin{defn}Let $p$ be a noncrossing partition of order $n$
corresponding to a permutation $\sigma$ which pairs points
$\rightbul{i}$ with $\leftbul{\sigma(i)}$.  Then
$\dual{p}$, which we call the
{\em dual} of~$p$,
is the noncrossing partition arising from the
noncrossing matching
which pairs
$\leftbul{i}$ with $\rightbul{(\sigma(i)-1)}$ (where
the subtraction is performed modulo $n$).
\end{defn}
In other words, $\dual{p}$ is the
partition that arises if $\rightbul{1}$ is relabeled $\leftbul{1}$,
$\leftbul{2}$ is relabeled $\rightbul{1}$, and so on cyclically
(pulling the first point $\leftbul{1}$ around to
to become the last point $\rightbul{n}$),
so that the equivalence relation is given by the complementary
connected components to those that defined $p$.
Referring to Figure \ref{hulls},
we see for example that the partition dual to
$\{ \{1,3,4,7\},\{2\},\{5,6\},\{8\} \}$ is
$\{ \{1,2\}, \{3\}, \{4,6\}, \{5\}, \{7,8\} \}$.
Note that the dual of the dual does not yield the same noncrossing
partition back again, but rather a relabeling of the original
by $1/n$ of a full rotation.

Non-crossing partitions have a natural partial order $\preceq$ by refinement, where
the ``least'' partition has $n$ blocks and the ``greatest'' has only
a single block containing all elements $1, \ldots, n$.  Duality
reverses this partial order: $p \preceq q$ if and only if
$\dual{q} \preceq \dual{p}$.  Refinement $\preceq$ gives $\Lambda_n$ the
structure of a lattice, where the join of two partitions
comes from taking the union of the two equivalence relations and then
extending by both transitivity and the noncrossing constraint---that is,
if $i < j < k < \ell$ and $i \sim k, j \sim \ell$, we require
$j \sim k$ as well.  The meet operation is similarly well-defined
(for example, as the anti-dual of the join of the duals).
The lattice $\Lambda_n$ is graded by $n$ minus the number of blocks
in the partition, and every maximal chain has the same length:
$q$ covers $p$ in the partial order if and only if
$p$ is obtained by splitting a block of $q$
into two blocks (in a noncrossing way), or in other words
if $p \preceq q$ and they differ by $1$ in the grading.
Figure \ref{lattice} illustrates the lattice of noncrossing
partitions of order 4.

\begin{figure} 
 \centering 
 \includegraphics[width=0.6\textwidth]{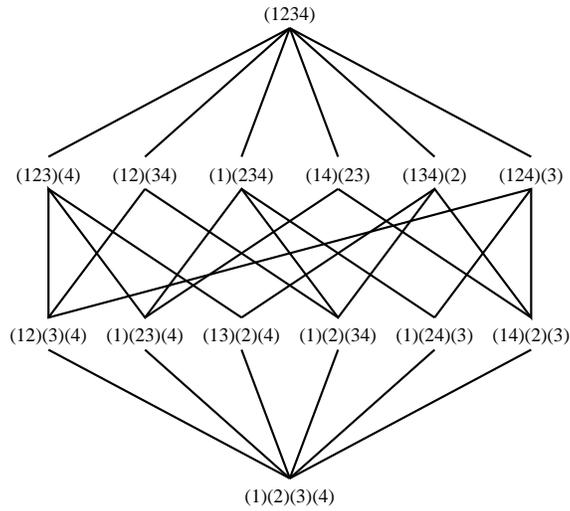} 
 \caption{The Hasse diagram of the lattice $\Lambda_4$.} 
 \label{lattice} 
\end{figure}
\end{subsection}

\begin{subsection}{The Cayley graph $\Gamma_n$ of transpositions in $S_n$}
One place where noncrossing partitions have turned up prominently
is in a particularly elegant solution to the word problem for
the braid group of order $n$.  Given any finite-type Artin group,
there are two canonical Garside structures, one arising from the
positive monoid on the standard generating set (the ``standard''
Garside structure), and one which is in many senses dual to it.
%
%
In the case of the braid group $B_n$, the
dual Garside monoid is generated by a set of $\binom{n}{2}$ conjugates
of the standard generators, whose image in the corresponding
Coxeter group $S_n$ is the set of transpositions (those permutations
that exchange two elements, leaving the others fixed).

We will consider permutations as vertices of the Cayley graph
$\Gamma_n$ of $S_n$ generated by all transpositions.
(Since these generators are involutions, we may take $\Gamma_n$
to be a simple, undirected graph.)
We write permutations in disjoint cycle notation,
optionally omitting fixed points.
We adopt the convention that permutations act on
the left, so that for example $(12)(23) = (123)$.
We give names to two distinguished vertices of $\Gamma_n$:
the identity permutation $e$, and the successor function
$s = (123\ldots n)$.  Given any two permutations $\sigma$
and $\tau$, we define the distance $d(\sigma,\tau)$ in
$\Gamma_n$ as the minimal length of a path from $\sigma$
to $\tau$.  By considering each edge of $\Gamma_n$ to
be a unit interval, we may extend $d$ to the entire graph.

The properties of $\Gamma_n$ as a Cayley graph
ensure that $d$ is symmetric and equivariant with respect
to left and right multiplication, so that in particular
$d(\sigma,\tau) = d(e,\sigma^{-1}\tau)$.  We also have the
following
\begin{prop}
If $\sigma \in S_n$ has $r$ orbits (including fixed points), then
$d(e,\sigma)$ = $n-r$.
\end{prop}
\begin{proof}
Let $\sigma \in S_n$ be written as a product of disjoint cycles
\[
(i_{11}i_{12}\ldots i_{1\ell_1})(i_{21}\ldots i_{2\ell_2})\ldots
(i_{r1}i_{r2}\ldots i_{r\ell_r}),
\]
including fixed points, so that $\ell_j$ may equal $1$ sometimes
and so that
$
\sum_{j=1}^r \ell_j = n.
$
On the one hand we have
\[
(i_{j1}\,i_{j2})(i_{j2}\,i_{j3})\ldots(i_{j\ell_j-1}\,i_{j\ell_j})
= (i_{j1}i_{j2}\ldots i_{j\ell_1}),
\]
so each cycle of length $\ell_j$ is a product of $\ell_j-1$ transpositions,
and
\[
d(e,\sigma) \leq \sum_{j=1}^r (\ell_j - 1) = n - r.
\]
But on the other hand multiplication of a permutation by a
transposition cannot reduce the number of cycles
by more than 1, so
$d(e,\sigma) \geq n - r$ and we have equality.
\end{proof}

The definition of duality for noncrossing partitions has a convenient
translation in terms of permutations:
If $\sigma$ is the permutation such that
$\rightbul{i}$ is paired with $\leftbul{\sigma(i)}$, then
the permutation such that
$\leftbul{i}$ is paired with $\rightbul{(\sigma(i) - 1)}$
is $\sigma^{-1} s$.  To revisit Figure \ref{hulls}
for example,
we have $(8)(65)(2)(7431)(12345678) = (12)(3)(46)(5)(78)$.
Applying the dual twice gives us
$(\sigma^{-1}s)^{-1}s = s^{-1}\sigma s$,
or in other words
conjugation by $1/n$ of a full rotation.
\end{subsection}
\end{section}

\begin{section}{results}
\begin{subsection}{The lattice $\Lambda_n$ of noncrossing partitions within
$\Gamma_n$}
Given a metric space $(X,d)$ and a pair of points $x,z \in X$, the {\em interval}
$\mathcal{I}_{x,z}$
from $x$ to $z$ is the union of all global geodesics between $x$ and $z$,
that is
\[ \mathcal{I}_{x,z} = \{y \in X : d(x,y) + d(y,z) = d(x,z) \}. \]

The following theorem is implicit in the standard developments
of the dual Garside structure on the braid group,
although it is not usually stated in terms of an interval
in the Cayley graph, nor is it usually proved
in terms of the genus of a constructed surface.
\begin{thm}
\label{single}
%
%
%
Let each vertex in the
Hasse diagram $\Lambda_n$ of the lattice of noncrossing partitions
be labeled with the permutation whose cycles are the blocks of
the partition, each cycle written in numerically increasing order.
Then this labeled Hasse diagram is the interval
$\mathcal{I}_{e,s}$
in the metric space $(\Gamma_n, d)$.
\end{thm}
\begin{proof}
The conventions we have chosen---that the noncrossing partition
corresponding to a
permutation
$\sigma$ connects~$\rightbul{i}$ to~$\leftbul{\sigma(i)}$,
and that the boundary points of a noncrossing partition
are labeled in counter-clockwise order around the disc---ensure
that the permutation $\sigma$ associated to a noncrossing
partition $p$ consists of a cycle for each block
of $p$, such that each cycle can be written in numerically increasing order.

If a permutation $\sigma$ does come from a non-crossing partition,
then the $n-1$ connected components of the disc divided by the
arcs of the noncrossing matching correspond to the cycles of
$\sigma$ and of the dual $\sigma^{-1}s$.  Each connected component
is a topological disc, and they are glued together along the arcs.
We imitate this construction in the case of an arbitrary
permutation, obtaining a disc exactly when the permutation
lies in $\Lambda_n$.

Let $\sigma$ be an arbitrary permutation, possibly crossing; we
extend the definition of the dual $\dual{\sigma}$
as $\sigma^{-1}s$.  We construct an oriented surface with boundary
as follows:
for each cycle $c = (c_1 c_2 \ldots c_k)$ of $\sigma$, create a topological disc
with $2k$ points marked
$\leftbul{c_1},\rightbul{c_1},\leftbul{c_2},\ldots,\rightbul{c_k}$ in
counterclockwise order, where of course $c_{i+1} = \sigma(c_i)$.
The $k$ edges from $\rightbul{c_i}$ to $\leftbul{\sigma(c_i)}$ we label
$c_i$ and orient counter-clockwise (that is,
in the same order as previously listed)
so that the union of discs for all cycles of $\sigma$ has,
for each label $1 \leq i \leq n$,
exactly one edge with the label $i$ going from the point marked
$\rightbul{i}$ to the point marked $\leftbul{\sigma(i)}$.
Similarly, for each cycle $c = (c_1 c_2 \ldots c_k)$
of $\dual{\sigma}$, we mark $2k$ points on the boundary of a disc
as $\rightbul{c_1},\leftbul{(c_1+1)},\rightbul{c_2},\ldots,\leftbul{(c_k+1)}$,
where $c_{i} = \sigma^{-1}s(c_{i-1})$ and thus $c_{i-1} + 1 = \sigma(c_{i})$.
The $k$ edges from
$\leftbul{\sigma(c_{i})}$ to $\rightbul{c_{i}}$
we label $c_i$ and orient clockwise (that is,
in the reverse order of that previously listed)
so that the union of discs for all cycles of $\dual{\sigma}$
also has a full complement of edges labeled $i$ and going from
$\rightbul{i}$ to $\leftbul{\sigma(i)}$.

If $\sigma$ represented a noncrossing partition to begin with, then
what we have just described in rather intricate detail is the
collection of connected components obtained by cutting a disc
along the arcs of the corresponding noncrossing matching---refer
for example to Figure \ref{hulls}.
For an example of what happens in the more general case,
refer to
Figure \ref{genus2} which illustrates the collection of
marked discs obtained from the permutation
$(15)(2436)$.

\begin{figure} 
 \centering 
 \includegraphics[width=0.6\textwidth]{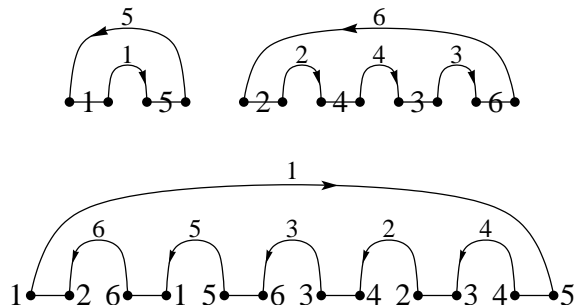} 
 \caption{The components of $\Sigma_{\sigma}$ for
$\sigma = (15)(2436)$ (top row from $\sigma$;
bottom row from $\dual{\sigma}$).}
 \label{genus2} 
\end{figure}

Identifying pairs of oriented edges labeled $i$ between the two
collections of discs, we obtain a connected oriented
surface $\Sigma_{\sigma}$
whose single boundary component has $2n$ vertices marked
$\leftbul{1}, \rightbul{1}, \ldots, \leftbul{n},\rightbul{n}$ in
order.  The permutation $\sigma$ represents a noncrossing partition
if and only if $\Sigma_{\sigma}$ is a disc, which is true
(for a connected compact surface with one boundary component)
if and only if the Euler characteristic of $\Sigma_{\sigma}$ is $1$.

We have $2n$ vertices for $\Sigma_{\sigma}$ and $3n$ edges ($2n$ on the boundary
and $n$ in the interior).  Each cycle of $\sigma$ contributes one
face, and since $d(e,\sigma) = n - k$ where $k$ is the number of cycles,
we have $n - d(e,\sigma)$ faces coming from $\sigma$.  Similarly, we have
$n - d(e,\sigma^{-1}{s}) = n - d(\sigma,s)$ faces coming from
$\dual{\sigma}$.  It follows that the Euler characteristic of $\Sigma_{\sigma}$
is $1$ exactly when \[d(e,\sigma) + d(\sigma,s) = n-1 = d(e,s)\]
and
thus that
the vertices of $\Lambda_n$ are exactly those
permutations which lie in $\mathcal{I}_{e,s}$.
Indeed, the genus of
$\Sigma_{\sigma}$ is
\[ \frac12\left[d(e,\sigma) + d(\sigma,s) - d(e,s)\right] \]
and measures exactly the degree to which $\sigma$ fails
to lie in the interval.

It remains only to show that the edges of $\Lambda_n$ are the same
as those of $\mathcal{I}_{e,s}$.  Every edge in $\mathcal{I}_{e,s}$
is part of a path of length $n-1$ from $e$ to $s$, and each step
of such a path unites two orbits of a permutation; this is a covering
relation in the refinement order $\preceq$ of noncrossing partitions.
Conversely, every covering relation in the order $\preceq$ joins two partition
blocks $(i_1 i_2 \ldots i_j)$ and $(i_{j+1} \ldots i_\ell)$ to
form $(i_1 i_2 \ldots i_\ell)$; this can be accomplished  by left
multiplication by the transposition $(i_1\,i_{j+1})$, which decreases
the distance to $s$ in $\Gamma_n$ and is thus an edge of
$\mathcal{I}_{e,s}$.
\end{proof}

\begin{rem}
Given any pair of permutations $\sigma,\tau \in \Gamma_n$ for which
$d(\sigma,\tau) = n-1$,
$\sigma^{-1}\tau$ consists of a single cycle (conjugate to $s$)
and $\mathcal{I}_{\sigma,\tau}$ is isomorphic to
$\mathcal{I}_{e,s}$ (since the generating set of $\Gamma_n$ is invariant
with respect to conjugation).
There are thus $n!(n-1)!$ isometrically embedded copies of
$\Lambda_n$ in $\Gamma_n$.
\end{rem}

Having defined $\Lambda_n$ in terms of the metric $d$ on $\Gamma_n$, it
is natural to define a metric $d_{\Lambda}(\sigma, \tau)$ as the length
of the shortest path from $\sigma$ to $\tau$ which is contained
entirely within $\Lambda_n$.  (In fact, as we shall shortly see,
$d_\Lambda = d$.)

\end{subsection}

\begin{subsection}{Meanders in $\Lambda_n$}
We now turn our attention to pairs of permutations $(\sigma,\tau)$
where each of the permutations $\sigma, \tau$ belongs to
$\Lambda_n = \mathcal{I}_{e,s}$.  When we construct the noncrossing
matching corresponding to $\sigma$ in the upper half-plane
and the (reflection of the) noncrossing matching corresponding
to $\tau$ in the lower half-plane, we obtain a
{\em system of meanders}, a collection of simple closed
curves with some number of components
$\Pi_0(\sigma,\tau)$.  The meanders are exactly those
pairs $(\sigma, \tau)$ for which $\Pi_0(\sigma,\tau) = 1$.

\begin{thm}
\label{pairs}
For any pair of permutations $\sigma, \tau \in \Lambda_n$,
the following are equivalent:
\begin{enumerate}
\item $\Pi_0(\sigma,\tau) = k$
\item $d(\sigma,\tau) = n - k$
\item $d_\Lambda(\sigma,\tau) = n - k$
\end{enumerate}
\end{thm}

\begin{proof}
Let $M$ be the system of meanders corresponding to the pair
$(\sigma, \tau)$, where $\sigma$ and $\tau$ are both permutations
belonging to the lattice $\Lambda_n$.  As usual, we label the
transverse crossings of $M$ as $\leftbul{1}, \rightbul{1},
\ldots, \leftbul{n}, \rightbul{n}$ in order.

As each odd point $\leftbul{i}$ is connected to an even point
by an arc in the upper half-plane, in order to count the
number of components $\Pi_0(\sigma,\tau)$ of the system of
meanders $M$ it suffices to find which components contain
which of the even points $\rightbul{i}$.  We define
a permutation of the integers $1,\ldots,n$ as follows:
starting at an even point $\rightbul{i}$, first follow
the arc connecting $\rightbul{i}$ in the lower half-plane
to $\leftbul{\tau(i)}$, then follow the arc in the upper
half-plane which connects $\leftbul{\tau(i)}$ to
$\rightbul{\sigma^{-1}\tau(i)}$.  Repeating this permutation
we will eventually traverse the entire connected component
of $M$ containing $\rightbul{i}$; it follows that the
components of $M$ are in one-to-one correspondence
with the cycles of the permutation $\sigma^{-1}\tau$.
We thus have
\[
d(\sigma,\tau) = d(e,\sigma^{-1}\tau) = n - k,
\]
where $k = \Pi_0(\sigma,\tau)$ is the number of connected
components of $M$.

Now to show the equivalence of $d$ and $d_\Lambda$.  Every
path in $\Lambda_n$ is a path in $\Gamma_n$, so clearly
$d_\Lambda(\sigma, \tau) \geq d(\sigma,\tau)$.  We need to
show that the distance between $\sigma$ and $\tau$
in $\Gamma_n$ can
be realized by a path lying entirely within $\Lambda_n$.

We proceed by induction.  For $n = 1$, we have of
necessity $\sigma = e$, $\tau = e$, and
$d_\Lambda(\sigma,\tau) = d(\sigma,\tau) = 0$.  For
$n > 1$, our aim will be to reduce to the following
case:  Suppose $\sigma(n) = n$ and $\tau(n) = n$
both.  Then $\sigma$ and $\tau$ can both be taken
to lie within $\Gamma_{n-1}$ and $\Lambda_{n-1}$,
and by induction a shortest path within $\Gamma_{n-1}$,
which is also a shortest path within $\Gamma_n$, can
be taken to lie within $\Lambda_{n-1} \subset \Lambda_n$.

To reduce to this case, consider first the action of the
duality map $(\sigma,\tau) \mapsto (\dual{\sigma},\dual{\tau})$.
This is an isometry of the Cayley graph $\Gamma_n$ which
exchanges the endpoints $e$ and $s$ of $\mathcal{I}_{e,s}$,
and so it is also an isometry of $\Lambda_n$, so
$d_\Lambda(\sigma,\tau) = d(\sigma,\tau)$ if and only if
$d_\Lambda(\dual{\sigma},\dual{\tau})
= d(\dual{\sigma},\dual{\tau})$.  The effect on $M$
of taking simultaneous duals on $\sigma$ and $\tau$
is to pull the first section of curve passing through
$\leftbul{1}$ around to become $\rightbul{n}$,
relabeling $\rightbul{1}$ as $\leftbul{1}$,
$\leftbul{2}$ as $\rightbul{1}$, and so forth.

Now, since the system of arcs in the lower half-plane
is planar and noncrossing, there must be an innermost
arc, an arc which encloses no other.  By applying
duality repeatedly, this innermost arc can be moved
to the first position connecting $\leftbul{1}$
with $\rightbul{1}$, and by applying duality twice
more we may assume without loss of generality
that $\tau(n) = n$.  If $\sigma(n) = n$ also, then
by induction we are done; otherwise we are halfway
there.

Assume then that $\tau(n) = n$ but $\sigma(n) = i \neq n$.
In $M$ we thus have a curve which descends through
point $\leftbul{n}$, bends around immediately to
pass through $\rightbul{n}$, and then continues on to
$\leftbul{i}$.
Let $\sigma^\prime$ be the product of permutations
$(i\, n)\,\sigma$; then the system
of meanders $M^\prime$ given by $(\sigma^{\prime},\tau)$ has
a circle through the points $\leftbul{n}$ and
$\rightbul{n}$, and the section of curve which used
to pass down through $\leftbul{n}$ instead
doubles back up to eventually pass
down through $\leftbul{i}$.
In particular,
\begin{itemize}
\item The permutation $\sigma^{\prime}$ corresponds
to the noncrossing partition obtained by splitting
out a singleton block $\{ n \}$ from the noncrossing
partition of $\sigma$, and is thus one step
away from $\sigma$ in $\Lambda_n$.
\item The system of meanders $M^\prime$ has exactly
one more component than $M$, so
$\Pi_0(\sigma^{\prime},\tau) = \Pi_0(\sigma,\tau) + 1$
\item and thus $d(\sigma^{\prime},\tau) = d(\sigma,\tau) -1$.
\item Since both $\tau(n) = n$ and $\sigma^\prime(n) = n$,
we can assume by induction that
$d_\Lambda(\sigma^\prime,\tau) = d(\sigma^\prime,\tau)$.
\end{itemize}
It follows that a path realizing the length $d(\sigma,\tau)$
can be taken entirely within $\Lambda_n$, the first step
of which is from $\sigma$ to $\sigma^\prime$.
\end{proof}


\begin{rem}
David Savitt has proven independently
\cite{S06}
that meanders correspond to
diameters of the Hasse diagram of the lattice of noncrossing partitions.
\end{rem}
\end{subsection}

\begin{subsection}{Open questions}
\begin{enumerate}
\item This gives a definition of ``meander'' for any finite Coxeter group $W$.
What can we say about these generalized meanders?
\item This suggests a generalization of the {\em meander determinant}.
Do the generalized determinants also factor as nicely?
\item Does this give any hope for the meander enumeration problem?
\end{enumerate}
\end{subsection}

\end{section}

\maketitle
\bibliographystyle{alpha}
\bibliography{ref}
\end{document}